\documentclass[12pt, oneside]{amsart}   	
\usepackage{geometry}
\usepackage{amsmath, amssymb}
\usepackage{epsfig}
\usepackage{epic,eepic}              		
\usepackage{graphicx}				
\newtheorem{thm}{Theorem}[section]
\newtheorem{prop}[thm]{Proposition}
\newtheorem{prob}[thm]{Problem}
\newtheorem{lem}[thm]{Lemma}

\newtheorem{conj}[thm]{Conjecture}
\newtheorem{rem}[thm]{Remark}

\theoremstyle{definition}
\newtheorem{defn}[thm]{Definition}

\numberwithin{equation}{section}

\title{On configurations where the Loomis-Whitney inequality is nearly sharp and applications to the Furstenberg set problem}
\author{Ruixiang Zhang}
\usepackage{fancyhdr}
\begin{document}

\maketitle
\begin{abstract}
In this paper, we consider the so-called ``Furstenberg set problem'' in high dimensions. First, following Wolff's work on the two dimensional real case, we provide ``reasonable'' upper bounds for the problem for $\mathbb{R}$ or $\mathbb{F}_p$. Next we study the ``critical'' case and improve the ``trivial'' exponent by $\Omega (\frac{1}{n^2})$ for $\mathbb{F}_p^n$. Our key tool to obtain this lower bound is a theorem about how things behave when the Loomis-Whitney inequality is nearly sharp, as it helps us to reduce the problem down to dimension two.
\end{abstract}
\section{introduction}
The well-known Kakeya conjecture, which plays a key role in many harmonic analysis problems, says:
\begin{conj}[Kakeya Conjecture]
A Kakeya set in $\mathbb{R}^n$ has full Hausdorff dimension.
\end{conj}

To understand this problem, many efforts have been put in and variants have been studied. One possible variant is to assume we don't have a ``whole line'' but instead ``some lower dimensional part of it'' in each direction and then to see what happens. Wolff\cite{wolff1999recent} and Tao\cite{taoedinburgh} raised the following problem:
\begin{prob}[Furstenberg set problem]\label{furstenbergprob}
Fix $0 < \beta  \leq 1$. If a compact set $S$ in $\mathbb{R}^n$ satisfies that for any direction $\omega \in S^{n-1}$, there is a line parallel to $\omega$ such that a $\beta$-dimensional subset of this line lies in $S$. Then what can we get as the best lower bound for the dimension of $S$?
\end{prob}

Wolff \cite{wolff1999recent} attributes Problem \ref{furstenbergprob}  to Furstenberg's work and it is likely that \cite{furstenberg1970intersection} inspired his formulation of the problem. He only asked the two dimensional case and obtained a lower bound $\max \{\beta + \frac{1}{2}, 2\beta\}$ and an upper bound $\frac{3}{2} \beta + \frac{1}{2}$. Tao\cite{taoedinburgh} then asked the question in arbitrary dimensions.

From Wolff's lower bound we see that when $n = 2$ in Problem \ref{furstenbergprob}, the case $\beta = \frac{1}{2}$ is ``critical'' and is of much interest. People studied the ``critical'' problem and its interesting connection to several other problems (see for example \cite{katz2001some}). They found evidences leading to the conjecture suggesting that improvements of the lower dimension bound $1$ might be possible.  This was achieved by Bourgain \cite{bourgain2003erdos} via a quantitative sum-product-type lemma about $\delta$-discretized sets. Once this was established he could use the reduction in \cite{katz2001some} to show the following theorem:

\begin{thm}[\cite{bourgain2003erdos}]\label{bourgainthm}
When $n = 2$ and $\beta = \frac{1}{2}$ in Problem \ref{furstenbergprob}, then $\dim S > 1 + \varepsilon$ for an absolute constant $\varepsilon > 0$.
\end{thm}

The high dimensional case $n > 2$ is much harder and not much study has been done. For the upper bound, one can try to generalize Wolff's result and we will indeed do this in the following section. We will prove an upper bound, which is also a reasonable conjectural lower bound, along Wolff's lines. Note that one may also look at the counterpart of Problem \ref{furstenbergprob} in finite fields. For those fields, we prove the same upper bounds with some additional assumption when the field is not of prime order. This, in particular, fully answers a question of Guth \cite{guth2012polynomial} modulo a constant. As Wolff already noticed, there can be better upper bounds when the finite field is not of prime order. For the prime order case we will again put our upper bound as the conjectural lower bound.

As for provable lower bounds over $\mathbb{R}$, one may use some standard methods like the Kakeya maximal inequality approach performed by Wolff \cite{wolff1999recent} to get some lower bounds. However, we don't yet know how to prove full Kakeya in high dimensions yet. This makes those bounds almost surely not good enough, let alone a ``beyond-Kakeya estimate'' as achieved in Theorem \ref{bourgainthm} in $2$ dimensions.

Finite fields have an advantage that Kakeya is easier. In fact, one major breakthrough in the past few years was the proof of the Kakeya conjecture for finite fields by Dvir\cite{dvir2009size}. He used the so-called polynomial method in the proof, making it surprisingly short and elegant.

A direct generalization of his polynomial argument yields a lower dimension bound $n \beta$ for the finite field analogue of Problem \ref{furstenbergprob}. On the other hand, by counting pairs of points lying on a same line one may obtain a trivial bound $\beta + \frac{n-1}{2}$. So the ``easy'' lower bound is $ \max \{ \beta + \frac{n-1}{2}, n\beta \}$ for any finite field. This bound agrees with the lower bound obtained by Wolff in dimension two. Also, from this we notice that $\beta = \frac{1}{2}$ can be viewed as the ``critical'' exponent for all $n$. Like discussions in \cite{bourgain2003erdos} and \cite{katz2001some}, we will be interested only in this critical setting for the lower bound, which is already hard (and is in a sense the ``hardest'' case because that the cases when $\beta = 0$ and $\beta = 1$ are easy, see the end of the next section). From now on we fix $\beta = \frac{1}{2}$ in this section.

A next interesting question is then to investigate whether we can go beyond the bound given by Kakeya to obtain things like Theorem \ref{bourgainthm}. Here one needs to be cautious: as Wolff\cite{wolff1999recent} already noticed, we couldn't go further in the $\mathbb{F}_{p^2}$ setting ($p$ is a prime in this paper). We will mention this phenomenon in the next section. Thus we will only seek for an improvement in $\mathbb{F}_p$. Here we can indeed go beyond the lower dimension bound $1$ in the case $n = 2$ by the Szemer\'{e}di-Trotter theorem in $\mathbb{F}_p$ proved in \cite{bourgain2004sum}. It's not surprising at all that some sum-product type ingredients come in, just like the case in \cite{katz2001some}\cite{bourgain2003erdos}.

For an arbitrary dimension $n$, we discover that we can indeed go beyond Kakeya in all dimensions in the $\mathbb{F}_p$ analogue of Problem \ref{furstenbergprob}. What we will show in this paper is:
\begin{thm}[Main Theorem]\label{mainthmFF}
Assume $p$ is a prime, $n \geq 2$. If a subset $S \subseteq \mathbb{F}_p^n$ satisfies that for any direction there is a line $l$ in this direction such that $|l \bigcap S| \geq p^{\frac{1}{2}}$, then $|S| \gtrsim p^{\frac{n}{2} + \Omega (\frac{1}{n^2})}$.
\end{thm}

Here $\Omega (\cdot)$ is the standard Landau asymptotic notation so that $f_1 = \Omega(f_2) (f_1, f_2 > 0)$ means that there exists a positive constant $\varepsilon$ s.t. $f_1 \geq \varepsilon f_2$. The main idea in our proof is to exploit the fact that the number of incidences is close to the trivial bound and do a projective transformation to put the points onto a ``grid''. We then make use that there are many essentially non-cohyperplanar lines passing through these points. The key observation is that we can analyze what happens when the equality of Loomis-Whitney is approached (which is our case) and then project everything to a $2$ dimensional subspace. There we could use Szemer\'{e}di-Trotter to get a nontrivial improvement of the exponent. Note that the conjectured upper bound that I suggest, as we shall see, is $p^{\frac{3n-1}{4}}$. So there is still a large gap in between.

In the following discussion, all the implied constants will depend solely on dimension $n$ (unless otherwise specified) and may differ even in a single chain of inequalities. By a refinement $T'$ of a finite set $T$, we mean that $T' \subseteq T$ and that there is a constant $\varepsilon > 0$ such that $|T'| \geq \varepsilon |T|$.

\section*{Acknowledgements}
I was supported by the mathematics department of Princeton University. I would like to thank Larry Guth who mentioned the importance of refinement to me when I tried to formulate Theorem \ref{LWeqthm}. His notes were also the ones bringing this problem to my attention. I would like to thank Terence Tao for mentioning to me the work of Furstenberg\cite{furstenberg1970intersection} that inspired Problem \ref{furstenbergprob}. Also, I would like to thank Jordan Ellenberg. The helpful discussion with him made me rethink about the upper bound.

\section{Upper bounds for $\mathbb{R}$ and $\mathbb{F}_q$}

In this section we prove ``reasonable'' upper bounds for the Furstenberg set problem in any dimension, both for $\mathbb{R}^n$ and for $\mathbb{F}_p^n$. Following Wolff's heuristic \cite{wolff1999recent} \cite{wolff2003lectures}, we also conjecture the upper bounds will be sharp for $\mathbb{R}^n$ since we can prove the same lower bound for a ``discrete model'' in $\mathbb{R}^n$ (see \cite{zhang2014polynomials}). The reader should notice that things are more complicated for general finite fields: For example, we have better configurations for, say, $\mathbb{F}_{p^2}$, which was already noticed \cite{taoedinburgh} \cite{wolff1999recent} \cite{wolff2003lectures}. Nevertheless, we do not have counterexamples for prime fields and I conjecture again that the upper bound should be sharp there.

We gradually demonstrate the above assertions. First look at the ``real'' problem in $\mathbb{R}^n$. Here when $n = 2$ Wolff proved an upper bound, which we state as a theorem below.

\begin{thm}[Wolff \cite{wolff1999recent}\cite{wolff2003lectures}]\label{Wolffexample}
Assume $n=2$. In Problem \ref{furstenbergprob}, the (Hausdorff) dimension of $S$ can be as small as $\frac{1}{2} + \frac{3}{2} \beta$.
\end{thm}

We use a direct generalization of the construction used by Wolff in the proof of Theorem \ref{Wolffexample} to show the following high dimensional generalization.

\begin{thm}\label{ourrealexample}
In Problem \ref{furstenbergprob}, the (Hausdorff) dimension of $S$ can be as small as $\frac{n-1}{2} + \frac{n+1}{2} \beta$.
\end{thm}

We need a version of Jarnik's theorem to prove Theorem \ref{ourrealexample}.

\begin{thm}\label{jarnikthm}
If $\{N_j\}$ is a sequence of positive integers and is increasing sufficiently rapidly, and $0< \beta \leq 1$, then the set
\begin{equation}
T = \{ x \in (\frac{1}{4}, \frac{3}{4}) : \text{ for any } j, \text{ there exist } p, q \in \mathbb{Z} \text{ s.t. } q \leq N_j^{\beta} \text{ and } |x - \frac{p}{q}| \leq N_j^{-2} \}
\end{equation}
has Hausdorff dimension $\beta$.
\end{thm}

The proof of Theorem \ref{jarnikthm} is immediate from the proof of another slightly different Jarnik-type theorem by Besicovitch \cite{besicovitch1934sets}. Thus we omit the proof.

\begin{proof}[Proof of Theorem \ref{ourrealexample}]
This proof is entirely parallel to Wolff's argument in dimension $2$ \cite{wolff1999recent} \cite{wolff2003lectures}.

Following Wolff's notations, denote a $G$-set in $\mathbb{R}^n$ to be a compact set $E \subseteq \mathbb{R}^n$ that is contained inside $[0,1] \times [-100, 100] \times [-100, 100] \times \cdots \times [-100, 100]$ such that for each $(m_2, m_3 , \ldots, m_n) \in [0, 1]^{n-1}$ there is a line segment contained in $E$ connecting $x_1 = 0$ to $x_1 = 1$ with ``slope'' $(m_2, m_3 , \ldots, m_n)$. In other words for each $(m_2, m_3 , \ldots, m_n) \in [0, 1]^{n-1}$ there exists a vector $(b_2, b_3 , \ldots, b_n)$ such that $(x, m_2 x + b_2, \ldots, m_n x + b_n) \in E, 0 \leq x \leq 1$.

With a line $l : (x, \mathbf{m} x + \mathbf {b})$, where $\mathbf{m} \neq 0$ and $\mathbf{b}$ are $n-1$ dimensional vectors, we associate a $\delta$-tube ${\mathbb{R}^n} \supseteq S_l^{\delta} = \{(x, \mathrm{y}): 0 \leq x \leq 1, |\mathbf{y}- (\mathbf{m} x + \mathbf{b})| \leq \delta\}$ around it.

By Theorem \ref{jarnikthm}, the set
\begin{equation}
T' = \{ t: \frac{1-t}{\sqrt{2}t} \in T\}
\end{equation}
also has dimension $\beta$.

For a large fixed $N$, we consider the family of all line segments $l_{j_2  j_3 \cdots j_n k_1 k_2 \cdots k_n}$ connecting a point $(0, \frac{j_2}{N}, \frac{j_3}{N}, \ldots, \frac{j_n}{N})$ to a point $(1, \frac{k_2}{N} \sqrt{2}, \frac{k_3}{N} \sqrt{2}, \ldots, \frac{k_n}{N} \sqrt{2})$ where $j_2 , \ldots, j_n, k_2, \ldots, k_n$ are integers between $0$ and $N$. Thus
\begin{equation}
l_{j_2  j_3 \cdots j_n k_1 k_2 \cdots k_n} = \{(x, \mathbf{\Phi}_{\mathbf{j}\mathbf{k}} (x))\}
\end{equation}
where the $n-1$ dimensional vectors $\mathbf{\Phi}_{\mathbf{j}\mathbf{k}} (x) = (1-x) \frac{\mathbf{j}}{N} + x \frac{\mathbf{k}}{N} \sqrt{2}$, $\mathbf{j} = (j_2, \ldots, j_n)$, $\mathbf{k} = (k_2, \ldots, k_n)$. For simplicity we also call $l_{j_2  j_3 \cdots j_n k_1 k_2 \cdots k_n}$ to be $l_{\mathbf{j}\mathbf{k}}$.

By e.g. Example 3.2 in p. 124 of \cite{kuipers2012uniform}, every vector in $[0, 1]^{n-1}$ has a distance $\lesssim N^{-2} (\log N)^2$  from the ``slope vector'' (i.e. $-\frac{\mathbf{j}}{N} + \frac{\mathbf{k}}{N} \sqrt{2}$) of one of the $l_{\mathbf{j}\mathbf{k}}$. Thus the set $G_N = \cup_{\mathbf{j}, \mathbf{k}} S_{l_{\mathbf{j}\mathbf{k}}}^{N^{-2} (\log N)^3}$ is a $G$-set.

Let
\begin{equation}
Q_N = \{t : \frac{1-t}{\sqrt{2} t} \text{ is a rational number } \frac{p}{q} \in (\frac{1}{4}, \frac{3}{4}) \text{ with denominator } q \leq N^{\beta}\}.
\end{equation}

For each $t \in Q_N$, let $S(t) = \{\mathbf{\Phi}_{\mathbf{j}\mathbf{k}} (t)\}_{j_2, \ldots, j_n, k_2, \ldots, k_n = 0}^{N}$. Then for each pair $\mathbf{j}$ and $\mathbf{k}$ we have that $\frac{\mathbf{\Phi}_{\mathbf{j}\mathbf{k}} (t)}{\sqrt{2} t} = \frac{p\mathbf{j} + q \mathbf{k}}{qN}$ is a rational vector of bounded size with a common denominator $qN$. Thus $|S(t)| \lesssim (qN)^{n-1} \lesssim N^{(n-1)(1+ \beta)}$. Hence $|\cup_{t \in Q_N} S(t)| \lesssim N^{n-1 + (n+1) \beta}$ and we deduce that: The set $H_N = \{(x, \mathbf{y}) \in G_N : |x -t| \leq \frac{1}{N^2} \text{ for  some } t \in Q_N\}$ is contained in the union of $\lesssim N^{n-1 + (n+1) \beta}$ disks of radius $N^{-2} (\log N)^3$.

Now we let a sequence $N_j$ increase rapidly. We recursively construct compact sets $F_j \subseteq \mathbb{R}^n$ such that $F_{j+1} \subseteq F_j$, each $F_j$ is a $G$-set and that the set $\widetilde{F_j} = \{(x, \mathbf{y}) \in F_j : x \in T'\}$ is contained inside a union of $\lesssim N_j^{n-1 + (n+1) \beta} \log N_j$ disks of radius $N_j^{-2} (\log N_j)^3$.

We proceed as follows: Let $F_0$ be any $G$-set. If $F_j$ has been constructed, it will have a form $F_j = \cup_{i=1}^M S_{l_i}^{\delta}$ for some $\delta$, where $l_i = \{(x, \mathbf{m}_i x + \mathbf{b}_i) : 0 \leq x \leq 1\}$ for some $\mathbf{m_i}, \mathbf{b}_i$ and every vector in $[0, 1]^{n-1}$ is $\leq \delta$ away from some $\mathbf{m}_i$. Let $A_i (x, \mathbf{y}) = (x, \mathbf{m}_i x + \delta \mathbf{y} + \mathbf{b}_i)$. We can choose $N_{j+1}$ to be sufficiently large and define
\begin{equation}
F_{j+1} = \cup_{i=1}^m A_i (G_{N_{j+1}}).
\end{equation}

It's plain that $F_{j+1} \subseteq F_j$ and it is not hard to check that $F_{j+1}$ is a $G$-set. The covering property is also essentially obvious from the covering properties of $H_N$ provided $N_{j+1}$ is large enough.

Finally take $F = \cap_{j=1}^{\infty} F_j$ and let $E = \{(x, \mathbf{y}) \in F : x \in T'\}$. By the covering properties of $\widetilde{F_j}$, $\dim E \leq \frac{n-1}{2} + \frac{n+1}{2} \beta$. Also,  $F$ itself is a $G$-set. Any line segment $l \subseteq F$ that connects the hyperplane $x_1 = 0$ and $x_1 = 1$ will satisfy $\dim (l \cap E) = \dim T' \geq \beta$. Take $S$ to be a finite union of rotated copies of $E$ and we complete the proof.
\end{proof}

\begin{rem}
When $\beta = 1$ this is the Kakeya setup. We can see that by this example we give nothing new for the Kakeya set problem.
\end{rem}

If we are looking at a ``discrete'' analogue of the Furstenberg set problem then the exponent $\frac{n-1}{2} + \frac{n+1}{2} \beta$ is indeed best possible. Namely, in \cite{zhang2014polynomials} we prove the following theorem.

\begin{thm}\label{maththeoremoffurstenberg}
Given any fixed number $C > 0$ and $0 \leq \beta \leq 1$. Assume that a set $L$ of lines in $\mathbb{R}^n$ satisfies that the direction set of $L$ is a $\frac{C}{N}$-dense subset of the unit sphere. If a point set $P$ has at least $N^{\beta}$ points lying on each line in $L$ then $|P| \gtrsim_{n , C} N^{\frac{n-1}{2} + \frac{n+1}{2} \beta}$. Moreover, the exponent $\frac{n-1}{2} + \frac{n+1}{2} \beta$ cannot be improved.
\end{thm}

This can be viewed as a ``discrete analogue'' of the Fursetnberg set problem and thus following Wolff's philosophy, the discrete result might suggest a plausible exponent to conjecture. Thus we conjecture:

\begin{conj}
In Problem \ref{furstenbergprob}, the (Hausdorff) dimension of $S$ can not be smaller than $\frac{n-1}{2} + \frac{n+1}{2} \beta$.
\end{conj}

For finite fields, we have some more subtle issues when the order of the field is not a prime. For example the following construction (which was already noticed in the two dimensional case by Wolff and is a straight generalization of the construction there) shows that we may have an upper bound beyond what Theorem \ref{ourrealexample} suggests.

\begin{prop}\label{ourexamplepsquare}
Assume $n \geq 2$ and $p$ is a prime. Then there exists a subset $S \subseteq \mathbb{F}_{p^2}^n$ satisfying that $|S| \lesssim p^n$, and that for any direction there is a line $l$ in this direction such that $|l \bigcap S| \geq p$.
\end{prop}

\begin{proof}
We will essentially use nothing more than a projective transform image of $\mathbb{F}_p^n \subseteq \mathbb{F}_{p^2}^n$. The readers who are convinced here can skip the proof.

Fix an irreducible quadratic polynomial $P(x)$ over $\mathbb{F}_p$ and identify $\mathbb{F}_{p^2}$ with $\mathbb{F}_p [x] / (P(x))$. Let $\mu = \overline{x} \in \mathbb{F}_{p^2}$. We construct a set
\begin{equation}
\mathbb{F}_{p^2}^n \supseteq E_1 = \{(\frac{b}{a \mu + b}, \frac{c_2}{a \mu + b}\mu, \ldots, \frac{c_{n}}{a \mu + b}\mu) : a, b, c_2, \ldots, c_{n} \in \mathbb{F}_p, a \text{ or } b \neq 0\}.
\end{equation}

The definition has a slight redundancy that the points in $E$ does not change if we multiply $a, b, c_1, \ldots, c_{n-1}$ simultaneously by any $0 \neq \lambda \in \mathbb{F}_p$. We deduce that $|E_1| \sim p^n$.

Next look at any line connecting the points $(0, u_2, \ldots, u_n)$ and $(1, v_2 \mu, \ldots, v_n \mu)$ in $E_1$. Here $u_j, v_j \in \mathbb{F}_p$. We have $p+1$ different points $(\frac{b}{a \mu + b}, \frac{au_2 + b v_2}{a \mu + b}\mu, \ldots, \frac{a u_n + b v_n}{a \mu + b}\mu)$ (where $(a, b) \in \mathbf{P}\mathbb{F}_p^1$) of $E_1$ on this line. Since $v_j \mu - u_j$ can be any element in $\mathbb{F}_{p^2}$, we deduce that for any vector $(1, \mathbf{y})$ where $\mathbf{y} \in \mathbb{F}_{p^2}^{n-1}$, we have a line which meets $E_1$ at $\geq p$ points and is parallel to the vector. It suffices to construct $E_2, \ldots, E_n$ similarly and take their union to be $S$.
\end{proof}

The exponent in Proposition \ref{ourexamplepsquare} is definitely worse than the analogous exponent in \ref{ourrealexample} when $\beta = \frac{1}{2}$. Also, working on general finite fields for a general range of $\beta$ can have some new technical difficulties. It's not our intention to get too involved here and we will stick to prime fields or certain exponents in the rest of this section. We get the counterpart of Theorem \ref{ourrealexample} as the following:

\begin{thm}\label{ourrealexampleprimefield}
Assume $n \geq 2$, $0 \leq \beta \leq 1$, $K$ is a given positive constant and $q$ is a  sufficiently large prime depending on $K$. Then there exists a subset $S \subseteq \mathbb{F}_q^n$ satisfying that $|S| \lesssim_K q^{\frac{n-1}{2} + \frac{n+1}{2} \beta}$, and that for any direction there is a line $l$ in this direction such that $|l \bigcap S| \geq K q^{\beta}$. When $\beta = 1$ we require $K \leq 1$. When $\beta = 0$ the above also holds when $q$ is a prime power (where the implied constant does not depend on the form of $q$).
\end{thm}

Note that this theorem fully answers Question 2 in \cite{guth2012polynomial} when ``$A$ is constant'' there (corresponding to the ``$\beta = 0$'' case here), since by counting we have this bound to be sharp, modulo a constant, when $\beta = 0 $.

The proof of Theorem \ref{ourrealexampleprimefield} has common features with Theorem \ref{ourrealexample}. We first prove a lemma.

\begin{lem}\label{denselem}
Given a positive number $K>1$ which will only be used to state the lemma when the $q$ below is not a prime. For a finite field $\mathbb{F}_q$, there exists a set $\Delta$ of elements such that:

(i) $|\Delta| \sim q^{\frac{1}{2}}$.

(ii)  There exists a $\mu \in \mathbb{F}_q$ such that $\mu \Delta - \Delta = \mathbb{F}_q$. Here $A - B = \{a-b: a \in A, b \in B\}$.

(iii) If $q = p$ is a prime, then for any given $p^{\frac{1}{2}} \lesssim r \leq p$, $L(r) = |\{x \in \mathbb{F}_p: x \neq - \mu, |x\Delta + \Delta| \leq r\}|$ is $\gtrsim \frac{r^2}{p}$. If $q$ is a prime power then there are $K$ different $x \in \mathbb{F}_q$ such that $|x \Delta + \Delta| \lesssim_K q^{\frac{1}{2}}$ (and thus automatically $x \neq -\mu$ when $q$ is large).

In this lemma all the implied constants are universal independent of $q$. When $q$ is a prime power the constants may depend on $K$.
\end{lem}

\begin{proof}
First we assume $q=p$ is a (sufficiently large) prime. Let $s = \lceil \sqrt{p} \rceil$ be the least integer $\geq \sqrt{p}$. We take $\Delta = \{ \overline{1}, \overline{2}, \ldots, \overline{s}\}$ and $\mu = \overline{s+1}$. Apparently (i) and (ii) are satisfied. (iii) is also not hard to verify if we notice that for positive integers $x, y \leq t < s-1$, $|\overline{\frac{x}{y}} \Delta + \Delta| = |\overline{x} \Delta + \overline{y} \Delta| \lesssim t q^{\frac{1}{2}}$. Of course for $0< x, y, x', y' < s-1$, $\overline{\frac{x}{y}} = \overline{\frac{x'}{y'}}$ if and only if $\frac{x}{y} = \frac{x'}{y'}$. Look at all the pairs $(x, y)$ such that $1 < x, y \leq t$ and that they produce different $\frac{x}{y}$. The total number of such pairs is $\gtrsim (\sum_{j=1}^t \frac{\phi(j)}{j})t \gtrsim t^2$.

Next we assume $q = p^m$ is a (large) prime power. Notice that we can not assume $p$ to be large. Assume $p^{\gamma - 1} < K \leq p^{\gamma}$. We have $\gamma = o(m)$ when $q$ is large. Take a fixed irreducible polynomial $P(x)$ of degree $m$ over $\mathbb{F}_p$ and from now on identify $\mathbb{F}_q$ with $\mathbb{F}_p [x] / (P(x))$. We have two different cases according to the parity of $m$.

If $m=2h$ is even, then we take
\begin{equation}
\mathbb{F}_q \supseteq \Delta = \{\overline{a_{h-1} x^{h-1}+ a_{h-2} x^{h-2} + \cdots + a_0} : a_{h-1}, \ldots, a_0 \in \mathbb{F}_p\}
\end{equation}
and  (i) is plain. Choose $\mu = \overline{x^h}$ we have (ii). To check (iii), we note that for $x \in X = \{\overline{b_{\gamma-1} x^{\gamma-1} + \ldots + b_0}: b_j \in \mathbb{F}_p\}$ we have $|x\Delta + \Delta| \lesssim p^{h + \gamma} \lesssim Kq^{\frac{1}{2}}$. Finally $|X| = p^{\gamma} \geq K$.

If $m = 2h+1$ we will combine the thoughts above. Take $\Delta (p)$ to be the $\Delta$ for $\mathbb{F}_p$ and a $\mu(p)$ such that $\mu(p) \Delta(p) - \Delta (p) = \mathbb{F}_p$. We take
\begin{equation}
\mathbb{F}_q \supseteq \Delta = \{\overline{a_h x^h + a_{h-1} x^{h-1}+ a_{h-1} x^{h-2} + \cdots + a_0} : a_h, a_{h-1}, \ldots, a_0 \in \mathbb{F}_p, a_0 \in \Delta(p)\}
\end{equation}
and (i) follows from the corresponding property of $\Delta(p)$ (and notice that the implied constant is independent of $h$). Take $\mu = x^h + \mu(p)$ and (ii) follows from the corresponding property of $\Delta(p)$ and $\mu(p)$. Checking (iii) is a little more involved. If $\gamma = 1$ (i.e. $p$ is large) then we can take elements in $\mathbb{F}_p$ to verify (iii) by the corresponding property of $\Delta (p)$. If $\gamma > 1$ then $p \leq K$. Thus $p = O_K (1 )$. So we can again use the set $X = \{\overline{b_{\gamma-1} x^{\gamma-1} + \ldots + b_0}: b_j \in \mathbb{F}_p\}$ to verify (iii).
\end{proof}

\begin{proof}[Proof of Theorem \ref{ourrealexampleprimefield}]
With Lemma \ref{denselem}, the proof is straightforward as we did in the proof of Theorem \ref{ourrealexample} and Proposition \ref{ourexamplepsquare}. From now on we assume that $q = p$ is a prime. The other case is essentially similar.

Take the $\Delta, \mu$ as in Lemma \ref{denselem}. Also by property (iii) there, we can find a set $X \subseteq \mathbb{F}_p$ such that $|x\Delta + \Delta| \lesssim_K p^{\frac{1+\beta}{2}}$ for any $x \in X$ and that $|X| = \lceil Kp^{\beta}\rceil$. Take
\begin{equation}
E_1 = \{(\frac{1}{a \mu + 1}, \frac{(au_2 + v_2)}{a\mu + 1} \mu, \ldots, \frac{(au_n + v_n)}{a\mu + 1} \mu) : a \in X, u_2, \ldots, u_n, v_2, \ldots, v_n \in \Delta\}.
\end{equation}

According to the size condition of $\Delta$ and $X$, we have $|E_1| \sim p^{\frac{n-1}{2} + \frac{n+1}{2} \beta}$. Moreover, we look at any line connecting two points $(0, u_2, \ldots, u_n)$ and $(1, v_2 \mu, \ldots, v_n \mu)$. It has $\geq Kp^{\beta}$ intersections with $E_1$ due to the choice of $X$. Since $v_j \mu - u_j$ can be any element in $\mathbb{F}_p$ we deduce that for any vector $(1, \mathbf{y})$ where $y \in \mathbb{F}_p^{n-1}$ , we have a line which meets $E_1$ at $\geq Kp^{\beta}$ points and is parallel to the vector. Now construct $E_2, \ldots, E_n$ similarly and take the union of the $n$ sets to be $S$.
\end{proof}

We conclude by conjecture this bound to be sharp for general prime fields.

\begin{conj}
Assume $n \geq 2$, $0 \leq \beta \leq 1$ and $p$ is a prime. $S \subseteq \mathbb{F}_p^n$ is a subset satisfying that for any direction there is a line $l$ in this direction such that $|l \bigcap S| \geq p^{\beta}$. Then $|S| \gtrsim q^{\frac{n-1}{2} + \frac{n+1}{2} \beta}$.
\end{conj}

 This is not even known for the dimension two case except when $\beta = 0$ or $1$. When $\beta = 0$ this is trivial by counting pairs of points on a same line. When $\beta = 1$ this is the finite field Kakeya which is true in general high dimensions and for general $\mathbb{F_q}$\cite{dvir2009size}. In general the two corresponding approaches together give a lower bound $q^{\max\{\beta + \frac{n-1}{2}, n \beta\}}$ (which is valid for all finite fields $\mathbb{F}_q$ and is not improvable for, e.g., $q = p^2$ and $\beta = \frac{1}{2}$). That is simply not enough here. A sum-product argument like \cite{bourgain2004sum} can improve the easy lower bound a little bit in dimension two. And we will generalize this gain to high dimensions in the following section.

\section{Proof of the main theorem}

In this section we prove Theorem \ref{mainthmFF}. Our main tool will be a theorem which illustrates how things are like when the Loomis-Whitney inequality is close to equality. In fact, this theorem implies that when that is the case, a positive proportion of the points have the ``smallest'' projection we can expect onto any fixed subspace. In applications, we will only use (\ref{conclusionoflemma2}) but we also write several other quite interesting statements that we could naturally obtain . We do this because we believe this theorem is interesting in its own right. For our convenience we introduce a definition before the theorem.

\begin{defn}
For a set $B$, a positive integer $n$ and a subset $T \subseteq B^n$, we denote the \emph{projection} of an element $x = (x_1, x_2, \ldots, x_n) \in T$ onto coordinate subspaces by the following:
\begin{equation}
Pr_{k_1, k_2, \ldots, k_t} (x) =  (x_{k_1}, x_{k_2}, \ldots, x_{k_t}), \text{ for any } t \leq n, 1\leq k_1 < k_2 < \cdots < k_t \leq n.
\end{equation}
Also we define
\begin{equation}
Pr_{k_1, k_2, \ldots, k_t} (T) =  \{Pr_{k_1, k_2, \ldots, k_t} (x): x \in T\}.
\end{equation}

We will often compose projections. To be convenient in this setting, if $j_1, \ldots, j_{t_1}$ is a subsequence of $k_1, \ldots, k_{t_2}$ we can define $Pr_{j_1, j_2, \ldots, j_{t_1}}$ (we use the same notation as the domain of this map will always be clear) from $Pr_{k_1, k_2, \ldots, k_{t_2}} (T)$ to $Pr_{j_1, j_2, \ldots, j_{t_1}} (T)$ as
\begin{equation}
Pr_{j_1, j_2, \ldots, j_{t_1}} (y) = \text{ the only element in } Pr_{j_1, j_2, \ldots, j_{t_1}} (Pr_{k_1, k_2, \ldots, k_{t_2}}^{-1} (y)),  \text{ for any } y \in Pr_{k_1, k_2, \ldots, k_{t_2}} (T).
\end{equation}
\end{defn}

\begin{thm}\label{LWeqthm}
Fix a set $B$ and positive integers $m \leq n$. Assume a subset $T \subset B^n$ satisfying the following ``size condition for projections'':
\begin{equation}
|Pr_{1, 2, \ldots, \widehat{k}, \ldots, n} (T)| \leq N,  \text{ for any } 1\leq k\leq n.
\end{equation}

Then there is a refinement $T'$ of $T$ such that
\begin{equation}\label{conclusionoflemma1}
|\{z \in Pr_{1, 2, \ldots, \widehat{j}, \ldots, n}(T') : Pr_{1, 2, \ldots, m} (z) = y\}| \gtrsim \frac{|T|^{n-m-1}}{N^{n-m-1}},  \text{ for any } y \in Pr_{1, 2, \ldots, m} (T') \text{ and }  m < j \leq n,
\end{equation}
\begin{equation}\label{conclusionoflemma3}
|Pr_{1, 2, \ldots, m}^{-1} (y) \bigcap T'| \gtrsim \frac{|T|^{n-m}}{N^{n-m}},  \text{ for any } y \in Pr_{1, 2, \ldots, m} (T'),
\end{equation}
\begin{equation}\label{conclusionoflemma2}
|Pr_{1, 2, \ldots, m} (T')| \lesssim \frac{N^{n-m}}{|T|^{n-m-1}},
\end{equation}
and
\begin{equation}\label{conclusionoflemma4}
|\{Pr_{1, 2, \ldots, \widehat{j}, \ldots, m} (T')\}| \lesssim \frac{N^{n-m+1}}{|T|^{n-m}},  \text{ for any } 1\leq j \leq m.
\end{equation}
\end{thm}

\begin{proof}
If $m = n$, we can just take $T' = T$.

If $m = n-1$, we can simply take a refinement $T' = \{x \in T: |Pr_{1, 2, \ldots, m}^{-1} (Pr_{1, 2, \ldots, m} (x))| \gtrsim \frac{|T|}{N}\}$. When the implied constant is small enough, it is easy to see this is indeed a refinement. Then it is trivial that (\ref{conclusionoflemma1}), (\ref{conclusionoflemma3}) and (\ref{conclusionoflemma2}) hold. The final part of the proof for $m \leq n-2$ (see below) will show that they together imply (\ref{conclusionoflemma4}), even in our case $m = n-1$.

Next we assume $m\leq n-2$. Suppose $Pr_{1, 2, \ldots, m} (T) = \{y_1, y_2, \ldots, y_s\}$. Define
\begin{equation}
a_j (y_k) = |\{z \in Pr_{1, 2, \ldots, \widehat{j}, \ldots, n} (T): Pr_{1, 2, \ldots, m} (z) = y_k\} |,  \text{ for any } m < j \leq n.
\end{equation}

Then we trivially have by assumption
\begin{equation}\label{conditionoflemma}
\sum_{k=1}^s a_j (y_k) \leq N,  \text{ for any } m < j \leq n.
\end{equation}

By Loomis-Whitney inequality,
\begin{equation}\label{LWlemmaestimate}
|Pr_{1, 2, \ldots, m}^{-1} (y_k)| \leq (\prod_{j= m+1}^n a_j (y_k))^{\frac{1}{n-m-1}}.
\end{equation}

Define
\begin{equation}\label{Ujdefine}
U_j = \{1 \leq k \leq s : (\prod_{l= m+1}^n a_l (y_k))^{\frac{1}{n-m-1}} \leq \frac{a_j (y_k) |T|}{100nN}\},  \text{ for any } m < j \leq n.
\end{equation}

Then for any $m < j \leq n$, by (\ref{conditionoflemma}), (\ref{LWlemmaestimate}) and (\ref{Ujdefine}),
\begin{equation}\label{forrefinementuse}
|Pr_{1, 2, \ldots, m}^{-1} (\{y_k : k \in U_j\})| = \sum_{k \in U_j} |Pr_{1, 2, \ldots, m}^{-1} (y_k)| \leq \sum_{k \in U_j} \frac{a_j (y_k) |T|}{100nN} \leq \frac{|T|}{100n}.
\end{equation}

We define $T_1= T - \bigcup_{j = m+1}^n Pr_{1, 2, \ldots, m}^{-1} (\{y_k : k \in U_j\})$. By (\ref{forrefinementuse}) it is a refinement of $T$. For any $x \in T_1$, assume $Pr_{1, 2, \ldots, m} (x) = y_k$, then $k$ is not in any $U_j , m< j \leq n$. Thus for any $m < j \leq n$,
\begin{equation}\label{conditionforTprime}
 (\prod_{l= m+1}^n a_l (y_k))^{\frac{1}{n-m-1}} \gtrsim \frac{a_j (y_k) |T|}{N}.
\end{equation}

If $a_{j_0} (y_k)$ is the smallest and $a_{j_1} (y_k)$ is the largest among all $a_j (y_k) (m<  j \leq n)$, without loss of generality we assume $j_0 \neq j_1$. Then (\ref{conditionforTprime}) implies
\begin{equation}
 (a_{j_0} (y_k))^{\frac{1}{n-m-1}} a_{j_1} (y_k)\gtrsim \frac{a_{j_1} (y_k) |T|}{N}.
\end{equation}

Therefore for any $m < j \leq n$ , $a_j (y_k) \gtrsim {(\frac{|T|}{N})}^{n-m-1}$, which implies that for any $j$, $|\{z \in Pr_{1, 2, \ldots, \widehat{j}, \ldots, n} (T') : Pr_{1, 2, \ldots, m} (z) = Pr_{1, 2, \ldots, m} (x)\}| \gtrsim {(\frac{|T|}{N})}^{n-m-1}$. This verifies (\ref{conclusionoflemma1}) for $T' = T_1$ and (\ref{conclusionoflemma2}) trivially follows for $T' = T_1$.

However, $T' = T_1$ may not satisfy (\ref{conclusionoflemma3}) and (\ref{conclusionoflemma4}). We refine it once more. Take a subset
\begin{equation}\label{chooseT2}
T_2 = \{ x \in T_1 : |Pr_{1, 2, \ldots, m}^{-1} (Pr_{1, 2, \ldots, m} (x))| \gtrsim \frac{|T|^{n-m}}{N^{n-m}}\}.
\end{equation}

From (\ref{conclusionoflemma2}) for $T' = T_1$ we deduce that in (\ref{chooseT2}) if we choose the implied constant sufficiently small then $T_2$ is a refinement of $T_1$. Obviously, $T' = T_2$ satisfies (\ref{conclusionoflemma1}), (\ref{conclusionoflemma3}) and (\ref{conclusionoflemma2}).

Note that for any $1 \leq j \leq m$ and any $w \in Pr_{1, 2, \ldots, \widehat{j}, \ldots, m} (T_2)$, consider $S_j (w, T_2) = \{z \in Pr_{1, 2, \ldots, \widehat{j}, \ldots, n} (T_2) : Pr_{1, 2, \ldots, \widehat{j}, \ldots, m} (z) =w\}$. We take any $x \in T_2$ such that $Pr_{1, 2, \ldots, \widehat{j}, \ldots, m} (x) = w$. Assume $Pr_{1, 2, \ldots, m} (x) = y$, then by (\ref{conclusionoflemma3}) we have that for $T' = T_2$, $|Pr_{1, 2, \ldots, m}^{-1} (y) \cap T_2| \gtrsim \frac{|T|^{n-m}}{N^{n-m}}$. Apparently, all elements of $|Pr_{1, 2, \ldots, m}^{-1} (y)|$ must have different image under $Pr_{m+1, m+2, \ldots, n}$. Thus they have different image under $Pr_{1, 2, \ldots, \widehat{j}, \ldots, n}$. But those images are all transformed by $Pr_{1, 2, \ldots, \widehat{j}, \ldots, m}$ to $Pr_{1, 2, \ldots, \widehat{j}, \ldots, m} (y) = Pr_{1, 2, \ldots, \widehat{j}, \ldots, m} (x) = w$. Thus $|S_j (w, T_2)| \gtrsim |Pr_{1, 2, \ldots, m}^{-1} (y)| \gtrsim \frac{|T|^{n-m}}{N^{n-m}}$. Since $|Pr_{1, 2, \ldots, \widehat{j}, \ldots, n} (T_2)| \leq N$, we deduce (\ref{conclusionoflemma4}) for $T' = T_2$. Thus it suffices to take $T' = T_2$.
\end{proof}

Theorem \ref{LWeqthm} enables us to ``descent'' from high dimensional setting to low dimensional setting. Indeed, in the following proof of Theorem \ref{mainthmFF} the key point is to reduce the problem to the dimension two case. There it will be further reduced to the Szemer\'{e}di-Trotter theorem in $\mathbb{F}_p$ proved in e.g. \cite{bourgain2004sum}.

\begin{proof}[Proof of Theorem \ref{mainthmFF}]
If $n=2$, without loss of generality we may assume $|S| \leq p^{\frac{3}{2}}$. Consider $p$ lines in $p$ different directions, each passing through $\geq \sqrt{p}$ points in $S$. By Szemer\'{e}di-Trotter on $\mathbb{F}_p$ \cite{bourgain2004sum}, we deduce
\begin{equation}
p^{\frac{3}{2}} \lesssim p + |S| + (p|S|)^{\frac{3}{4} -\varepsilon}.
\end{equation}
which implies the conclusion.

Next we assume $n \geq 3$. Suppose $|S| = M$. Since the trivial bound for $M$ (which follows from counting pairs of points on a same line) is already $p^{\frac{n}{2}}$, for our conclusion we only need to seek a little bit of gain in power. By assumption there is a set $L$ of lines such that for each direction there is exactly one $l \in L$ pointing towards this direction and for each $l \in L$ there is a subset $S_l \subseteq S \bigcap l$ such that $|S_l| \sim p^{\frac{1}{2}}$ (we can take $|S_l| = \lceil p^{\frac{1}{2}} \rceil$). Take
\begin{equation}\label{defnofS1}
S_1 = \{x\in S : |\{l\in L : x \in S_l\}| \lesssim \frac{M}{\sqrt{p}}\}.
\end{equation}

Note that $|\{(x, l): x \in S_l\}| \sim p^{n-\frac{1}{2}}$. We can choose the implied constant in (\ref{defnofS1}) sufficiently large such that $|S - S_1| \leq \frac{p^n}{100M}$. Hence we have
\begin{equation}
|\{(x, x', l): x, x' \in S_1 \bigcap S_l\}| \gtrsim p^n.
\end{equation}

Take $S_2 \subseteq S_1$ such that
\begin{equation}\label{defnofS2}
S_2 = \{ x \in S_1: |\{l \in L: x \in S_l\}| \gtrsim \frac{p^{n-\frac{1}{2}}}{M} \text{ and } |\{(x', l): x, x' \in S_1 \bigcap S_l\}| \gtrsim \frac{p^n}{M}\}.
\end{equation}

Apparently if we choose the implied constants sufficiently small (in fact we only need to care about the latter condition since the former follows automatically from it), we can still have
\begin{equation}\label{S2mainineq}
|\{(x, x', l): x, x' \in S_2 \bigcap S_l\}| \gtrsim p^n.
\end{equation}

Heuristically, what we want to do next is to choose $n$ points such that a lot of points lie on some same $S_l$ with any one of the $n$ chosen points. Then we may perform a projective transform like in \cite{katz2001some} to map the $n$ chosen points to infinity points of $n$ coordinate directions, respectively. Then since each of these $n$ points lie on $\lesssim \frac{M}{\sqrt{p}}$ (close to $p^{\frac{n-1}{2}}$ when $M$ is close to $p^{\frac{n}{2}}$) lines in $L$, the image of a lot of points in $S$ will have to be confined in a ``small'' range determined by their projections onto coordinate hyperplanes. One finds this is very close to the equality of Loomis-Whitney and Theorem \ref{LWeqthm} will help us to control the projection of them to a 2-dimensional coordinate subspace. Finally we use the $2$ dimensional result (essentially Szemer\'{e}di-Trotter) to get a non-trivial gain of power.

If we want to map $n$ points to infinity points in $n$ coordinate directions, respectively, then we don't want them to lie on a $n-2$ dimensional subspace. Motivated by this, we classify the points in $S_2$ as the following. We call a point $x \in S_2$ \emph{hyperplanar} if and only if among all $l$ satisfying $x \in S_l$, we have $\geq \delta \frac{p^{n-\frac{1}{2}}}{M}$ lines that lie on a hyperplane. We call other points in $S_2$ \emph{non-hyperplanar}. Here $\delta = \delta (p, n)$ is a positive number \emph{depending on $p$ and $n$} we will soon choose.

For any non-hyperplanar point $x \in S_2$, we consider $F_x = \{ (x_1, x_2, \ldots, x_n) \in S_1^n :  \text{ for any } 1\leq k \leq n, \text{ there exists an } l_k \in L \text{ s.t. } x, x_k \in S_1 \bigcap S_{l_k}\}$. Suppose $r = r(x; x_1 , \ldots, x_n)$ is the dimension of the affine subspace $V = V(x; x_1, \ldots, x_n)$ spanned by the lines $xx_1, xx_2, \ldots, xx_n$.

Fix any $r_0 < n$, we count the tuples with $r  = r_0$. For any such tuple $(x_1, \ldots, x_n) \in F_x$, there are $r_0$ points among $x_1 , \ldots, x_n$ such that $V(x; x_1, \ldots, x_n)$ is spanned by the $r_0$ lines joining each one of them and $x$ (respectively). Without loss of generality we only care about those tuples with $V(x; x_1, \ldots, x_n)$ spanned by $xx_1, xx_2, \ldots, xx_{r_0}$ since the result of counting will be differ by no more than a constant $C (n)$. If we fix $x_1, x_2, \ldots, x_{r_0}$ in this tuple then all other points $x_{r_0 +1}, \ldots, x_n$ are on the $r_0$-dimensional affine subspace determined by $x_1. x_2, \ldots, x_{r_0}$ and $x$ (thus in any hyperplane containing it). By the assumption that $x$ is non-hyperplanar, we have $\leq \delta \frac{p^n}{M}$ choice of each one of $x_{r_0 +1}, \ldots, x_n$. Hence the number of tuples with $r = r_0 < n $ is $\lesssim \delta^{n-r_0} M^{2r_0 -n}p^{n(n-r_0)}$.

By (\ref{defnofS2}), $|F_x| \gtrsim \frac{p^{n^2}}{M^{n}}$. Elementary calculation shows that when $M \leq \frac{1}{100000} p^{\frac{n}{2} + \frac{1}{100n^2}}$ (we can make this assumption or we are done) we can choose $\delta$ being a small number times $p^{-\frac{1}{n}}$ such that $\frac{p^{n^2}}{M^{n}} \geq X\sum_{1 \leq r_0 \leq n-1} \delta^{n-r_0} M^{2r_0 -n}p^{n(n-r_0)}$ for any fixed $X > 0$. Thus for an arbitrary $x \in S_2$ being non-hyperplanar,
\begin{equation}\label{nonhyperplaneq}
|\{(x_1, \ldots, x_n) \in F_x : xx_1, xx_2, \ldots, xx_n \text{ span the whole space} \}| \geq \frac{|F_x|}{2}.
\end{equation}

Let $S_3 = \{x \in S_2 : x \text{ is hyperplanar}\} \subseteq S_2$, $S_4 = S_2 - S_3$ is the set of non-hyperplanar points. We distinguish two cases according to whether $|S_3| \gtrsim \frac{p^n}{M}$. Here the implied constant is so small that in the ``$\lesssim$'' case by (\ref{S2mainineq}) we can still have
\begin{equation}\label{S4mainineq}
|\{(x, x', l): x, x' \in S_4 \bigcap S_l\}| \gtrsim p^n.
\end{equation}

Before we proceed we remark that we expect that being non-hyperplanar is a generic case since we are facing a set containing a ``half-dimension'' subset of a line in each direction. So perhaps we can have some easy gain of power in the first case when the set is ``highly hyperplanar''. We will soon see it is indeed the case here.

If $|S_3| \gtrsim \frac{p^n}{M}$, we fix an affine hyperplane $V_x$ for each $x \in S_3$ such that there are $\geq \delta \frac{p^{n-\frac{1}{2}}}{M}$ lines in $L \bigcap \{\text{lines in }V_x\}$ passing through $x$. Summing over $x \in S_3$, we get $\gtrsim \delta \frac{p^{2n- \frac{1}{2}}}{M^2}$ lines in $L$ in total (note: each fixed line can be counted multiple times here, and in fact $O(\sqrt{p})$ times). We choose an $x_1 \in S_3$ arbitrarily. Then on $V_{x_1}$ there are $\geq \delta \frac{p^n}{M}$ points in $S$. There are $\leq p^{n-2}$ directions on this affine hyperplane and (when counted with multiplicity) those directions are counted $O(p^{n- \frac{3}{2}})$ times in total above (we use the fact that lines in $L$ have different directions). Thus as long as $\delta \frac{p^{2n- \frac{1}{2}}}{M^2} \gtrsim p^{n- \frac{3}{2}}$ with the implied constant sufficiently large (as we may assume), we can choose another point $x_2 \in S_3$ such that there are $\gtrsim \delta \frac{p^{2n- \frac{1}{2}}}{M^3}$ lines $\in L$ that are on $V_{x_2}$ and transverse to $V_{x_1}$. Those lines must thus intersect $V_{x_1}$ at $\leq 1$ point and hence give rise to $\gtrsim \delta \frac{p^{2n}}{M^3}$ new points in $S$ in total.

We can find $x_3, x_4, \ldots$ successfully in this way until one of the following happens: (a) We have chosen $\sim \sqrt{p}$ points so that any new line transversal to all previous affine hyperplanes cannot always provide $\sim \sqrt{p}$ new points. (b) We have chosen enough points such that the number of total directions on all chosen $V_{x_k}$ times $\sqrt{p}$ is $\gtrsim \delta \frac{p^{2n- \frac{1}{2}}}{M^2}$, this is equivalent to having chosen $\gtrsim \delta \frac{p^{n+1}}{M^2}$ points. Thus we could terminate after choosing $\sim \min\{\sqrt{p}, \delta \frac{p^{n+1}}{M^2}\}$ points $x_k$ to ensure that each time after the first we obtain $\gtrsim \delta \frac{p^{2n}}{M^3}$ new points in $S$ and finally we get $\sim \delta \frac{p^{2n}}{M^3} \min\{\sqrt{p}, \delta \frac{p^{n+1}}{M^2}\}$ points $\in S$. Thus
\begin{equation}
M \gtrsim \delta \frac{p^{2n}}{M^3} \min\{\sqrt{p}, \delta \frac{p^{n+1}}{M^2}\}
\end{equation}
and $M \gtrsim p^{\frac{n}{2} + \Omega (1)}$. We are done.

From now on we assume $|S_3| \lesssim \frac{p^n}{M}$ with a suitable implied constant such that (\ref{S4mainineq}) holds. For each $x \in S_4$ denote $G_x = \{x' \in S_1 : \exists l \in L \text{ s.t. } x, x' \in S_l\}$. By (\ref{S4mainineq}),
\begin{equation}
\sum_{x \in S_4} |G_x| \geq p^n
\end{equation}

Note that $|F_x| = |G_x|^n$. By H\"{o}lder, we deduce
\begin{equation}
\sum_{x \in S_4} |F_x| \geq \frac{p^{n^2}}{M^{n-1}}.
\end{equation}

By (\ref{nonhyperplaneq}) and pigeonhole, we deduce that there exists a $(x_1, \ldots, x_n) \in S_1^n$ such that there are $\gtrsim \frac{p^{n^2}}{M^{2n-1}}$ different $x \in S_4$ satisfying $F_x \supseteq (x_1, \ldots, x_n)$ and $xx_1, xx_2 ,\ldots, xx_n$ do not lie on a hyperplane.

We now perform a projective transformation $\mathcal{P}$ that maps $x_1, x_2, \ldots, x_n$ to points at infinity in the $n$ coordinate directions, respectively. This is possible since the existence of any above $x$ shows that $x_1, \ldots, x_n$ are not in any $n-2$ dimensional subspace simultaneously. For any $x$ satisfying above conditions, assume $y =\mathcal{P} (x)$. Then by assumption $y$ is not on the infinity hyperplane.  Since $x_1, \ldots, x_n \in S_1$, we can find some $K \lesssim \frac{M}{\sqrt{p}}$ and lines $h_{jk} ( 1\leq k \leq K)$ parallel to the $j$th coordinate direction for each $1 \leq j \leq n$, such that each $y$ above lies on one $h_{1 k_1}$, one $h_{2 k_2}$, $\ldots$, and one $h_{n k_n}$. Also, there are $\gtrsim \frac{p^{n^2}}{M^{2n-1}}$ different $y$ in total. We label them as $y_1, y_2, \ldots, y_r (r \gtrsim \frac{p^{n^2}}{M^{2n-1}})$.

 Apply (\ref{conclusionoflemma2}) in Theorem \ref{LWeqthm} $n \choose 2$ times, we may assume that some $r \lesssim r'  \leq  r$ such that for any fixed $1 \leq j_1 < j_2 \leq n$, the $j_1$th and $j_2$th coordinate (as a pair) of $y_1 , \ldots, y_{r'}$ has $\lesssim \frac{(\frac{M}{\sqrt{p}})^{n-2}}{(\frac{p^{n^2}}{M^{2n-1}})^{n-3}} = \frac{M^{2n^2 - 6n + 1}}{p^{n^3 - 3n^2 + \frac{1}{2} n - 1}}$ possibilities.

Consider the image of $L$ under $\mathcal{P}$. $\mathcal{P} (L)$ is a union of $p^{n-1}$ lines. Each $y_j (1 \leq j \leq r')$ has $\gtrsim \frac{p^{n- \frac{1}{2}}}{M}$ \emph{strong} incidences with $\mathcal{P} (L)$ by (\ref{defnofS2}) (in this statement we use the \emph{strong} sense of incidence: we say $\mathcal{P} (x)$ incident with $\mathcal{P} (l)$ \emph{strongly}  if and only if $x \in S_l$). We have $\gtrsim \frac{p^{n^2 + n -\frac{1}{2}}}{M^{2n}}$ strong incidences in total. Thus by refining we could easily see that there are $\gtrsim \frac{p^{n^2 + n - 1}}{M^{2n}}$ lines  in $\mathcal{P} (L)$, each passing through $\gtrsim \frac{p^{n^2 + \frac{1}{2}}}{M^{2n}}$ points among $y_j (1 \leq j \leq r')$.

For each line $\mathcal{P} (l)$ in $\mathcal{P} (L)$, obviously there exists $1 \leq j_1 < j_2 \leq n$ such that the projection of $\mathcal{P} (l)$ onto the $j_1, j_2$-coordinate subspace ($2$ dimensional subspace spanned by the $j_1$th and $j_2$th coordinate vector) is not a point. Without loss of generality, we may assume  that there are $\gtrsim \frac{p^{n^2 + n - 1}}{M^{2n}}$ lines  in $\mathcal{P} (L)$, each passing through $\gtrsim \frac{p^{n^2 + \frac{1}{2}}}{M^{2n}}$ points among $y_j (1 \leq j \leq r')$ and none of them is vertical to the $1, 2$-coordinate subspace.

Now use an orthogonal projection $\mathcal{Q}$ to project everything (strictly speaking, $\mathcal{P}(S)$ and $\mathcal{P} (L)$) onto the $1, 2$-coordinate subspace. For any selected line $\mathcal{P} (l)$, there are $\leq p^{n-2}$ selected lines whose projected image is $\mathcal{Q} (\mathcal{P} (l))$ (if a line can be projected to $\mathcal{Q} (\mathcal{P} (l))$ by $\mathcal{Q}$, then it must lie on a fixed hyperplane and so does its preimage before doing transform $\mathcal{P}$. But there are $\leq p^{n-2}$ such lines in $L$ that can lie on this hyperplane). Therefore considering the image under $\mathcal{Q}$ we obtain $\gtrsim \frac{p^{n^2  + 1}}{M^{2n}}$ lines and $\lesssim \frac{M^{2n^2 - 6n + 1}}{p^{n^3 - 3n^2 + \frac{1}{2} n - 1}}$ points in $\mathbb{F}_p^2$ such that each line passes through $\gtrsim \frac{p^{n^2 + \frac{1}{2}}}{M^{2n}}$ points. By Szemer\'{e}di-Trotter in finite fields\cite{bourgain2004sum},
\begin{equation}
\frac{p^{n^2  + 1}}{M^{2n}} \cdot \frac{p^{n^2 + \frac{1}{2}}}{M^{2n}} \lesssim \frac{p^{n^2  + 1}}{M^{2n}} + \frac{M^{2n^2 - 6n + 1}}{p^{n^3 - 3n^2 + \frac{1}{2} n - 1}} + (\frac{p^{n^2  + 1}}{M^{2n}} \cdot \frac{M^{2n^2 - 6n + 1}}{p^{n^3 - 3n^2 + \frac{1}{2} n - 1}})^{\frac{3}{4} -\varepsilon}.
\end{equation}

By a simple calculation we deduce $M \gtrsim p^{\frac{n}{2} + \Omega(\frac{1}{n^2})}$.
\end{proof}

\begin{rem}
Will the projective-transform-and-projection method in this paper (but certainly has its root in earlier works like \cite{wolff1999recent} and \cite{katz2001some}) help us to go beyond Kakeya in $\mathbb{R}^n (n  \geq 3)$? We don't know yet. One major difficulty, for instance, is a good analogue of $S_3$ here that enables us to do a good projective transform when $|S_3|$ is small and still have an argument to gain in the exponent when $|S_3|$ is large. We leave this challenging task to interested readers.
\end{rem}

\bibliographystyle{amsalpha}
\bibliography{stageref}

\providecommand{\bysame}{\leavevmode\hbox to3em{\hrulefill}\thinspace}
\providecommand{\MR}{\relax\ifhmode\unskip\space\fi MR }
\providecommand{\MRhref}[2]{%
  \href{http://www.ams.org/mathscinet-getitem?mr=#1}{#2}
}
\providecommand{\href}[2]{#2}
\begin{thebibliography}{BKT04}

\bibitem[Bes34]{besicovitch1934sets}
AS~Besicovitch, \emph{Sets of fractional dimensions (iv): on rational
  approximation to real numbers}, Journal of the London Mathematical Society
  \textbf{1} (1934), no.~2, 126--131.

\bibitem[BKT04]{bourgain2004sum}
Jean Bourgain, Nets Katz, and Terence Tao, \emph{A sum-product estimate in
  finite fields, and applications}, Geometric \& Functional Analysis GAFA
  \textbf{14} (2004), no.~1, 27--57.

\bibitem[Bou03]{bourgain2003erdos}
Jean Bourgain, \emph{On the erdos-volkmann and katz-tao ring conjectures},
  Geometric and Functional Analysis \textbf{13} (2003), no.~2, 334--365.

\bibitem[Dvi09]{dvir2009size}
Zeev Dvir, \emph{On the size of kakeya sets in finite fields}, Journal of the
  American Mathematical Society \textbf{22} (2009), no.~4, 1093--1097.

\bibitem[Fur70]{furstenberg1970intersection}
Harry Furstenberg, \emph{Problems in analysis (sympos. salomon bochner,
  princeton univ., princeton, n.j., 1969), pp. 41--59}, Princeton Univ. Press,
  Princeton, N.J. (1970).

\bibitem[Gut]{guth2012polynomial}
Larry Guth, \emph{The polynomial method, fall 2012, project list},
  http://math.mit.edu/\~{} lguth/PolyMethod/projlist.pdf.

\bibitem[KN12]{kuipers2012uniform}
Lauwerens Kuipers and Harald Niederreiter, \emph{Uniform distribution of
  sequences}, Courier Dover Publications, 2012.

\bibitem[KT01]{katz2001some}
Nets Katz and Terence Tao, \emph{Some connections between falconer¡¯s distance
  set conjecture and sets of furstenburg type}, New York J. Math \textbf{7}
  (2001), 149--187.

\bibitem[Tao]{taoedinburgh}
Terence Tao, \emph{Edinburgh lecture notes on the kakeya problem}.

\bibitem[TV06]{tao2006additive}
Terence Tao and Van~H Vu, \emph{Additive combinatorics}, vol.~13, Cambridge
  University Press, 2006.

\bibitem[Wol99]{wolff1999recent}
Thomas Wolff, \emph{Recent work connected with the kakeya problem, prospects in
  mathematics (princeton, nj, 1996), 129--162}, Amer. Math. Soc., Providence,
  RI \textbf{12} (1999).

\bibitem[Wol03]{wolff2003lectures}
Thomas~H Wolff, \emph{Lectures on harmonic analysis}, vol.~29, AMS Bookstore,
  2003.

\bibitem[Zha14]{zhang2014polynomials}
Ruixiang Zhang, \emph{Polynomials with dense zero sets and discrete models of
  the kakeya conjecture and the furstenberg set problem}, arXiv preprint
  arXiv:1403.1352 (2014).

\end{thebibliography}

Department of Mathematics, Princeton University, Princeton, NJ 08540

ruixiang@math.princeton.edu

\end{document}